\documentclass[a4paper]{mathscan}

\def \ff {\frac}

\def \la {\langle}
\def \ra {\rangle}

\def \R  {\hbox{\bf R}}

\newtheorem{theorem}{Theorem}

\newtheorem{remark}{Remark}

\newtheorem{cor}{Corollary}

\begin{document}
\title{Jensen's operator inequality and its converses}

\begin{abstract}
We give a general formulation of Jensen's operator inequality for unital fields
of positive linear mappings, and we consider different types of converse inequalities.
\end{abstract}

\author{Frank Hansen}
\address{Department of Economics, Copenhagen University,
Studiestraede 6, 1455 Copenhagen K,
DENMARK}\email{frank.hansen@econ.ku.dk}
\author{Josip Pe\v {c}ari\'{c}}
\address{Faculty of Textile Technology, University of Zagreb, Pierottijeva 6, 10000
Zagreb, CROATIA}\email{pecaric@hazu.hr}
\author{Ivan Peri\'{c}}
\address{Faculty of Chemical Engineering and Technology, University of Zagreb,
Maruli\' cev trg 19, 10000, Zagreb, CROATIA}\email{iperic@pbf.hr}
\maketitle

\section{Introduction}

Let $I$ be a real interval of any type. A continuous function $f:I\to\R$
is said to be operator convex if
\begin{equation}\label{OP1}
f\left(\lambda x+(1-\lambda )y\right)\leq \lambda f(x)+(1-\lambda
)f(y)\end{equation} holds for each $\lambda\in [0,1]$ and every pair of
self-adjoint operators $x$ and $y$ (acting) on an infinite dimensional
Hilbert space $H$ with spectra in $I$ (the ordering is defined by setting
$ x\le y $ if $ y-x $ is positive semi-definite).

Let $ f $ be an operator convex function defined on an interval $ I. $
Ch.~Davis \cite{D} proved\footnote{There is small typo in the proof. Davis
states that $ \phi $ by Stinespring's theorem can be written on the form $
\phi(x)=P\rho(x)P $ where $ \rho $ is a $ * $-homomorphism to $ B(H) $ and
$ P $ is a projection acting on $ H. $ In fact, $ H $ may be embedded in a
Hilbert space $ K $ on which $ \rho $ and $ P $ acts. The theorem then
follows by the calculation $ f(\phi(x))=f(P\rho(x)P)\le
Pf(\rho(x))P=P\rho(f(x))P=\phi(f(x)), $ where the pinching inequality,
proved by Davis in the same paper, is applied.}
 a Schwarz inequality
\begin{equation}\label{OP3}
f\left(\phi (x)\right)\leq \phi\left( f(x)\right),
\end{equation}
where $ \phi\colon\mathcal A\to B(H)$ is a unital completely
positive linear map from a $ C^* $-algebra $ \mathcal A $ to
linear operators on a Hilbert space $ H, $ and $x$ is a
self-adjoint element in $ \mathcal A$ with spectrum in $I.$
Subsequently M.~D. Choi \cite{Ch} noted that it is enough to
assume that $ \phi $ is unital and positive. In fact, the
restriction of $ \phi $ to the commutative $ C^*$-algebra
generated by $ x $ is automatically completely positive by a
theorem of Stinespring.

F. Hansen and G. K. Pedersen \cite{HP1} proved a Jensen type inequality
\begin{equation}\label{OP2}
f\left(\sum_{i=1}^{n}a_{i}^{*}x_ia_i\right)\leq
\sum_{i=1}^{n}a_{i}^{*}f(x_i)a_i
\end{equation}
for operator convex functions $f$ defined on an interval $ I=[0,\alpha) $
(with $ \alpha\le\infty $ and $ f(0)\le 0) $ and self-adjoint operators $
x_1,\dots,x_n $ with spectra in $I$ assuming that
$\sum_{i=1}^{n}a_{i}^{*}a_i={\bf 1}.$ The restriction on the interval and
the requirement $ f(0)\le0 $ was subsequently removed by B. Mond and J.
Pe\v cari\' c in \cite{MP}, cf.~also \cite{HP2}.

The inequality (\ref{OP2}) is in fact just a reformulation of
(\ref{OP3}) although this was not noticed at the time. It is nevertheless
important to note that the proof given in \cite{HP1} and thus the statement of the theorem, when
restricted to $ n\times n $ matrices, holds for the much richer class of $
2n\times 2n $ matrix convex functions.
Hansen and Pedersen used (\ref{OP2})
to obtain elementary operations on functions, which leave invariant the class of operator monotone functions.
These results then served as the basis for a new proof of L{\"o}wner's theorem
applying convexity theory and Krein-Milman's theorem.

Finally B. Mond and J. Pe\v cari\' c \cite{MP1} proved the inequality
\begin{equation}\label{OP4}
f\left(\sum_{i=1}^{n}w_{i}\phi_{i}(x_i)\right)\leq
\sum_{i=1}^{n}w_{i}\phi_i(f(x_i))
\end{equation}
for operator convex functions $ f $ defined on an interval $ I, $ where
$\phi_i :B(H)\to B(K)$ are unital positive linear maps, $x_1,\dots,x_n $
are self-adjoint operators with spectra in $I$ and $w_1,\dots,w_n $
are non-negative real numbers with sum one.

The aim of this paper is to find an inequality which contains
(\ref{OP2}), (\ref{OP3}) and (\ref{OP4}) as special cases. Since the
inequality in (\ref{OP4}) was the motivating step for obtaining converses
of Jensen's inequality using the so called Mond-Pe\v cari\' c method, we
also give some results pertaining to converse inequalities in the new
formulation.

\section{The main result}

\noindent {\bf Continuous Fields of Operators}

Let $ T $ be a locally compact Hausdorff space and let $ \mathcal A
$ be a $ C^* $-algebra. We
say that a field $ (x_t)_{t\in T} $ of operators in $ \mathcal A $
is continuous if the function $ t\to x_t $ is norm continuous on $
T. $ If in addition $\mu$ is a Radon measure on $T$ and the function
$t\to \|x_t\| $ is integrable, then we can form the Bochner integral
$\int_T x_t\,d\mu(t)$, which is the unique element in $ \mathcal A $
such that
\[
\varphi \left(\int_T x_t\, d\mu (t)\right) = \int_T \varphi (x_t)\, d\mu
(t)
\] for every linear functional $ \varphi $ in the norm dual $ \mathcal
A^*, $ cf. \cite[Section 4.1]{HP2}.

Assume furthermore that there is a field $ (\phi_t)_{t\in T} $ of
positive linear mappings $ \phi_t:\mathcal A\to\mathcal B $ from $
\mathcal A $ to another $ C^* $-algebra $ \mathcal B. $
We say that such a field is continuous if
the function $ t\to\phi_t(x) $ is continuous for every $
x\in\mathcal A. $ If the $ C^* $-algebras are unital and
the field $t\to\phi_t(\mathbf 1) $
is integrable with integral $\mathbf 1,$ we say that $(\phi_t)_{t\in
T}$ is {\it unital.}

\begin{theorem}\label{IntJens} Let $f:I\to\R $ be an operator convex
function defined on an interval $ I, $ and let $ \mathcal A $ and
$ \mathcal B $ be unital $ C^* $-algebras. If $ (\phi_t)_{t\in T} $ is a
unital field of positive linear mappings $ \phi_t:\mathcal A\to
\mathcal B $ defined on a locally compact Hausdorff space $ T $
with a bounded Radon measure $ \mu, $ then the inequality
\begin{equation}\label{OP6}
f\left(\int_T\phi_t(x_t)\,d\mu(t)\right)\leq \int_{T}\phi_t(
f(x_t))\,d\mu(t)
\end{equation}
holds for every bounded continuous field $ (x_t)_{t\in T}$ of self-adjoint
elements in $\mathcal A$ with spectra contained in $I.$
\end{theorem}

\begin{proof}
We first note that the function $ t\to\phi_t(x_t)\in\mathcal B $
is continuous and bounded, hence integrable with respect to the
bounded Radon measure $ \mu. $ We may organize the set $
\text{\it CB}(T,\mathcal A) $ of bounded continuous functions on $ T $ with
values in $ \mathcal A $ as a normed involutive algebra by
applying the point-wise operations and setting
\[
\|(y_t)_{t\in T}\|=\sup_{t\in T}\|y_t\|\qquad (y_t)_{t\in T}\in
\text{\it CB}(T,\mathcal A),
\]
and it is not difficult to verify that the norm is already complete and
satisfy the $ C^* $-identity. In fact, this is a standard construction in
$ C^* $-algebra theory. It follows that $ f((x_t)_{t\in T})=(f(x_t))_{t\in
T} $. We then consider the mapping
\[
\pi\colon \text{\it CB}(T,\mathcal A)\to  {\mathcal B}
\]
defined by setting
\[
\pi\left((x_t)_{t\in T}\right)=\int_T \phi_t(x_t)\, d\mu(t),
\]
and note that it is a unital positive linear map. Setting $ x=(x_t)_{t\in
T}\in \text{\it CB}(T,\mathcal A), $ we use inequality (\ref{OP3}) to obtain
\[
\begin{array}{rl}
\displaystyle f\left(\pi\left((x_t)_{t\in T}\right)\right)&=f(\pi(x))\le\pi(f(x))=\pi\left(f\bigl((x_t)_{t\in
T}\bigr)\right)\\[2ex]
&=\pi\left(\bigl(f(x_t)\bigr)_{t\in T}\right),
\end{array}
\]
which is the statement of the theorem.
\end{proof}
To illustrate various techniques in proving Jensen's operator
inequality, we give two proofs of Theorem \ref{IntJens} in the discrete
case $T=\{1,\ldots ,n\}$.
\begin{proof}
By using the continuity of $f$
and uniform approximation of self-adjoint operators by simple operators,
we may assume that $ x_i $ has a spectral resolution on the form
\[
x_i=\sum_{j\in J_i}\,t_{i,j}e_{i,j}\qquad i=1,\dots,n,
\]
where each $J_i$ is a finite set and $ \sum_{j\in J_i} e_{i,j}=\mathbf 1. $ We then have
\begin{eqnarray*}
\lefteqn{f\left(\sum_{i=1}^{n}\phi_i
(x_i)\right)=f\left(\sum_{i=1}^{n}\phi_i \left(\sum_{j\in
J_i}t_{i,j}e_{i,j}\right)\right)} \nonumber \\
& & =f\left(\sum_{i=1}^{n}\sum_{j\in
J_i}t_{i,j}\phi_{i}(e_{i,j})\right)=f\left(\sum_{i=1}^{n}\sum_{j\in
J_i}\sqrt{\phi_{i}(e_{i,j})}\;t_{i,j}\;\sqrt{\phi_{i}(e_{i,j})}\right)
\nonumber \\
& & \leq \sum_{i=1}^{n}\sum_{j\in
J_i}\sqrt{\phi_{i}(e_{i,j})}\;f(t_{i,j})\;\sqrt{\phi_{i}(e_{i,j})}=
\sum_{i=1}^{n}\phi_{i}\left(\sum_{j\in J_i}f(t_{i,j})e_{i,j}\right) \nonumber
\\
& & =\sum_{i=1}^{n}\phi_i \left(f(x_i)\right).
\end{eqnarray*}
{\it The second proof:} We use an the idea from \cite{FK}, confer
also \cite{MP}. If $f$ is operator convex in $I=[0,1)$ and
$f(0)\leq 0$, then there is a connection $\sigma$ such that
$-f(t)=t\,\sigma\, (1-t)$. We use the following two properties of a
connection, cf. \cite{AV,FK}.
\begin{enumerate}\itemsep=1ex
\item  $\phi\left(a\,\sigma\, b\right)\leq \phi (a)\,\sigma\,\phi (b)$\\
for positive linear maps $\phi$ and positive operators $a$ and $b.$

\item $\sum_{i=1}^{n}a_i\,\sigma\, b_i\leq
\left(\sum_{i=1}^n a_i\right)\,\sigma\,\left(\sum_{i=1}^n b_i\right)$\\
for positive $n$-tuples $(a_1,...,a_n)$ and $(b_1,...,b_n).$\hfill (subadditivity)

\end{enumerate}
We then obtain
\begin{eqnarray*}
\lefteqn{-\sum_{i=1}^{n}\phi_{i}(f(x_i))=\sum_{i=1}^{n}\phi_{i}
\left(x_i\,\sigma\,({\bf 1}-x_i)\right)} \nonumber \\
& & \leq \sum_{i=1}^{n}\phi_i(x_i)\,\sigma\,\phi_i({\bf 1}-x_i)\leq
\left(\sum_{i=1}^n\phi_i
(x_i)\right)\,\sigma\,\left(\sum_{i=1}^n\phi_i ({\bf 1}-x_i)\right)
\nonumber \\
& & =\left(\sum_{i=1}^n\phi_i(x_i)\right)\,\sigma\,\left({\bf
1}-\sum_{i=1}^n\phi_i(x_i)\right)=-f\left(\sum_{i=1}^n\phi_i(x_i)\right).
\end{eqnarray*}
Consider now an arbitrary operator convex function $f$ defined on $[0,1)$.
The function $\tilde{f}(x)=f(x)-f(0)$ satisfies $ \tilde f(0)=0 $ hence
\begin{equation}
f\left(\sum_{i=1}^{n}\phi_i (x_i)\right)-f(0)\mathbf 1\leq \sum_{i=1}^{n}\phi_i
(f\left(x_i\right))-f(0)\sum_{i=1}^{n}\phi_i ({\bf 1})
\end{equation}
from which the statement follows. We finally obtain (\ref{OP6}) in this setting for a
function $ f $ defined on an arbitrary interval $[\alpha ,\beta)$ by
considering the function $g(x)=f\left((\beta -\alpha )x+\alpha\right)$ on the interval $[0,1)$.
\end{proof}

Inequality (\ref{OP6}) is obviously a generalization of
the inequalities (\ref{OP2}), (\ref{OP3}) and (\ref{OP4}).

\section{Converses of Jensen's inequality}

The following theorem should be compared with Theorem 2.3 in
\cite{FHPS}.
For a function $f:[m,M]\to\R$ we use the standard notation:
$$\alpha_f=\ff{f(M)-f(m)}{M-m}\quad\text{and}\quad\beta_f=\ff{Mf(m)-mf(M)}{M-m}.$$

We will apply functions $ F(t,s) $ of two real variables to
operators. In simple cases, when it makes sense, we may just replace
numbers with operators. If for example $ F(t,s)=s^{-1}t^{1/2}s^{-1},
$ we may set $ F(x,y)=y^{-1} x^{1/2} y^{-1}. $ Otherwise we may use
the functional calculus on tensor products, see for example
\cite{kn:hansen:2003:1}.

\begin{theorem}\label{Rev1} Let $(x_t)_{t\in T}$ be a bounded continuous field
of self-adjoint elements in a unital
$C^{*}$-algebra ${\mathcal A}$ with spectra in $[m,M]$ defined on a locally
compact Hausdorff space $T$ equipped with a Radon measure $\mu,$ and let
$(\phi_t)_{t\in T}$ be a unital field of positive linear maps $\phi_t :{\mathcal
A}\to {\mathcal B}$ from ${\mathcal A}$ to another unital
$C^{*}$-algebra ${\mathcal B}$. Let $f,g:[m,M]\to\R$ and
$F:U\times V\to\R$ be functions such that
$f\left([m,M]\right)\subset U,$ $g\left([m,M]\right)\subset V$ and
$F$ is bounded. If $F$ is operator monotone in the first
variable and $f$ is convex in the interval $[m,M]$, then
\begin{equation}\label{OP7}
\begin{array}{l}
\displaystyle F\left[\int_{T}\phi_{t}\left(f(x_t)\right)d\mu(t),
g\left(\int_{T}\phi_{t}(x_t)d\mu (t)\right)\right]\\[3ex]
\leq\displaystyle
\sup_{m\leq z\leq M}F\left[\alpha_{f}z+\beta_{f},g(z)\right]{\bf 1}.
\end{array}
\end{equation}
In the dual case (when $f$ is concave) the opposite inequality
holds in (\ref{OP7}) with $\inf$ instead of $\sup$.
\end{theorem}
\begin{proof} For convex $f$ the inequality $f(z)\leq \alpha_fz+\beta_f$
holds for every $z\in [m,M]$. Thus, by using functional calculus,
$f(x_t)\leq \alpha_f x_t+\beta_f {\bf 1}$ for every $t\in T$.
Applying the positive linear maps $\phi_t$ and integrating, we obtain
\[
\int_{T}\phi_t\left(f(x_t)\right)d\mu (t)\leq \alpha_f\int_{
T}\phi_t(x_t)d\mu (t)+\beta_f {\bf 1}.
\]
Now, using
operator monotonicity of $F(\cdot ,v)$, we obtain
\begin{eqnarray*}
\begin{array}{l}
\lefteqn{F\left[\int_{T}\phi_{t}\left(f(x_t)\right)d\mu
(t),g\left(\int_{T}\phi_{t}(x_t)d\mu (t)\right)\right]}\\[3ex]
\displaystyle\leq F\left[\alpha_f\int_{T}\phi_{t}(x_t)d\mu (t)+\beta_f {\bf 1},
g\left(\int_{T}\phi_{t}(x_t)d\mu (t)\right)\right]\\[3ex]
\displaystyle\leq \sup_{m\leq z\leq M}F\left[\alpha_{f}z+\beta_{f},g(z)\right]{\bf 1}.
\end{array}
\end{eqnarray*}
\end{proof}

Numerous applications of the previous theorem can be given (see
\cite{FHPS}). We give  generalizations of some results from
\cite{T}.

\begin{theorem}\label{Rev0} Let $(A_t)_{t\in T}$ be a continuous field
 of positive operators on a Hilbert space $H$
defined on a locally compact Hausdorff space $T$ equipped with a Radon measure
$\mu. $ We assume the spectra are in $[m,M]$ for some $0<m<M.$ Let furthermore
$(x_t)_{t\in T}$ be a
continuous field of vectors in $H$ such that
$\int_{T}\|x_t\|^2d\mu (t)=1.$  Then for any $\lambda \geq 0$, $p\geq 1$ and $q\geq 1$ we have
\begin{equation}\label{Rev2} \left(\int_{T}\la A_{t}^p
x_t,x_t\ra d\mu(t) \right)^{1/q}-\lambda \int_{T}\la
A_tx_t,x_t\ra d\mu (t)\leq C(\lambda,m,M,p,q),
\end{equation}
where the constant
\begin{equation}
\begin{array}{l}
C(\lambda,m,M,p,q)\\[2ex]
=\left\{\begin{array}{ll}
M\left(M^{\ff{p}{q}-1}-\lambda\right), & 0<\lambda
\leq\ff{\alpha_p}{q}M^{p\left(\ff{1}{q}-1\right)} \\
\ff{q-1}{q}\left(\ff{q}{\alpha_p}\lambda\right)^{\ff{1}{1-q}}+\ff{\beta_p}{\alpha_p}\lambda,
 &\ff{\alpha_p}{q} M^{p\left(\ff{1}{q}-1\right)}\le \lambda
 \le \ff{\alpha_p}{q} m^{p\left(\ff{1}{q}-1\right)} \\
 m\left(m^{\ff{p}{q}-1}-\lambda\right), &\ff{\alpha_p}{q} m^{p\left(\ff{1}{q}-1\right)}\leq \lambda
\end{array}\right.
\end{array}
\end{equation}
and $ \alpha_p $ and $ \beta_p $ are the constants $ \alpha_f $ and $ \beta_f $ associated with the function
$ f(z)=z^p. $
\end{theorem}
\begin{proof} Applying Theorem \ref{Rev1} for the functions
\[
f(z)=z^p,\quad F(u,v)=u^{1/q}-\lambda v,
\]
and unital fields of positive linear maps $ \phi_t\colon B(H)\to\mathbf C $
defined by setting $ \phi_t (A)=\la A x_t,x_t\ra $ for $ t\in T, $
the problem is reduced to determine
$\sup_{m\leq z\leq M}H(z)$ where
$H(z)=(\alpha_p z+\beta_p)^{1/q}-\lambda z$.
\end{proof}

The following Corollary is a generalization of Theorem 5 in
\cite{T}. The $r-$geometric mean $A\#_{r}B$ introduced by F. Kubo
and T. Ando in \cite{KA} is defined by setting
$$A\#_{r}B=A^{1/2}\left(A^{-1/2}BA^{-1/2}\right)^{r}A^{1/2}$$
for positive invertible operators $A$ and $B.$

\begin{cor}\label{Rev3} Let $(A_t)_{t\in T}$ and $(B_{t})_{t\in T}$ be
continuous fields of
positive invertible operators on a Hilbert space $H$ defined on a locally compact Hausdorff space $T$
equipped with a Radon measure $ \mu $ such that
\[
m_1{\bf 1}\leq A_t\leq M_1{\bf 1}\quad\text{and}\quad m_2{\bf 1}\leq B_t\leq
M_2{\bf 1}
\]
for all $ t\in T $ for some $ 0< m_1< M_1 $ and $ 0< m_2< M_2. $
Then for any $\lambda\ge 0$, $s\geq 1$, $p\ge 1$ and any continuous
field $(x_t)_{t\in T}$ of vectors in $H$ such that
$\int_{T}\|x_t\|^2 d\mu (t)=1$ we have
\begin{equation}\label{Rev33}
\begin{array}{l}
\displaystyle\left(\int_{T}\la A_{t}^p x_t,x_t\ra d\mu
(t)\right)^{1/p}\left(\int_{T}\la B_{t}^q x_t,x_t\ra d\mu(t)\right)^{1/q}\\[3ex]
\displaystyle\hskip 12em -\lambda\int_{T}\la
B_{t}^q\#_{1/s}A_{t}^px_t,x_t\ra d\mu (t)\\[3ex]
\displaystyle\leq C\left(\lambda,\ff{m_{1}^{p/s}}{M_{2}^{q/s}}\,,
\ff{M_{1}^{p/s}}{m_{2}^{q/s}}\,,s,p\right)M_{2}^{q}
 \end{array}
\end{equation}
where the constant $C$ is defined in Theorem \ref{Rev0} and $1/p+1/q=1.$
\end{cor}

\begin{proof} By using Theorem \ref{Rev0} we obtain for any $\lambda
\ge 0$, for any continuous field  $(C_t)_{t\in T}$ of positive
operators with $m{\bf 1}\leq C_t\leq M{\bf 1}$
and a square integrable continuous field $(y_t)_{t\in T}$ of vectors in $H$ the inequality
\begin{eqnarray}\label{Rev4}
\begin{array}{l}
\displaystyle\left(\int_{T}\la
 C_{t}^{s}y_t,y_t\ra\, d\mu (t) \right)^{1/p}\left(\int_{T}\la
y_t,y_t\ra\, d\mu (t)\right)^{1/q}\\[3ex]
\displaystyle\hskip 11em -\lambda\int_{T}\la
C_ty_t,y_t\ra\, d\mu (t)\\[3ex]
\displaystyle\leq C(\lambda,m,M,s,p)\int_{T}\la y_t,y_t\ra\, d\mu (t).
\end{array}
\end{eqnarray}
Set now $C_t=\left(B_{t}^{-q/2}A_{t}^pB_{t}^{-q/2}\right)^{1/s}$  and
$y_t=B_{t}^{q/2}x_t$ for $t\in T$ in (\ref{Rev4}) and observe that
$$\ff{m_{1}^{p/s}}{M_{2}^{q/s}}{\bf 1}\leq
\left(B_{t}^{-q/2}A_{t}^pB_{t}^{-q/2}\right)^{1/s}\leq
\ff{M_{1}^{p/s}}{m_{2}^{q/s}}{\bf 1}.$$
By using the definition of the
$1/s-$geometric mean and rearranging (\ref{Rev4}) we obtain
\[
\begin{array}{l}
\displaystyle\left(\int_{T}\la A_{t}^p x_t,x_t\ra\, d\mu (t)\right)^{1/p}
 \left(\int_{T}\la B_{t}^q x_t,x_t\ra\,  d \mu(t)\right)^{1/q}\\[3ex]
\hskip 12em\displaystyle-\lambda\int_{T}\la
B_{t}^q\#_{1/s}A_{t}^px_t,x_t\ra\, d\mu (t) \nonumber \\[2ex]
\displaystyle\leq  C\left(\lambda
,\ff{m_{1}^{p/s}}{M_{2}^{q/s}},\ff{M_{1}^{p/s}}{m_{2}^{q/s}},s,p\right)\int_{T}\la
B_{t}^qx_t,x_t\ra\, d\mu (t)\\[4ex]
\displaystyle  \leq  C\left(\lambda
,\ff{m_{1}^{p/s}}{M_{2}^{q/s}},\ff{M_{1}^{p/s}}{m_{2}^{q/s}},s,p\right)
M_{2}^{q}
\end{array}
\]
which gives (\ref{Rev33}).
\end{proof}

In the present context we may obtain results of the Li-Mathias
type by using Theorem \ref{Rev1} and the following result which is a simple consequence
of Theorem \ref{IntJens}.

\begin{theorem}\label{LM} Let $(x_t)_{t\in T}$ be a bounded continuous field
of self-adjoint elements in a unital $C^{*}$-algebra ${\mathcal A}$
defined on a locally compact Hausdorff space $T$ equipped with a Radon measure $ \mu. $
We assume the spectra are in $[m,M].$ Let furthermore
$(\phi_t)_{t\in T}$ be a unital field of positive linear maps $\phi_t :{\mathcal
A}\to {\mathcal B}$ from ${\mathcal A}$ to another unital
$C^{*}$-algebra ${\mathcal B}$. Let $f,g:[m,M]\to\R$ and
$F:U\times V\to\R$ be functions such that
$f\left([m,M]\right)\subset U,$ $g\left([m,M]\right)\subset V$
and $F$ is bounded. If $F$ is operator monotone in the first
variable and $f$ is operator convex in the interval $[m,M]$, then
\begin{equation}
\begin{array}{l}
\displaystyle F\left[\int_{T}\phi_{t}\left(f(x_t)\right)\,d\mu
(t),g\left(\int_{T}\phi_{t}(x_t)\,d\mu (t)\right)\right]\\[4ex]
\displaystyle\geq
\inf_{m\leq z\leq M}F\left[f(z),g(z)\right]{\bf 1}.
\end{array}
\end{equation}
In the dual case (when $f$ is operator concave) the opposite
inequality holds  with $\sup$ instead of $\inf$.
\end{theorem}

We also give generalizations of some results from \cite{DK}.

\begin{theorem}\label{DK1} Let $f$ be a convex function
on $[0,\infty )$ and let $\|\cdot\|$ be a normalized unitarily
invariant norm on $B(H)$ for some finite dimensional Hilbert space $H.$
Let $(\phi_t)_{t\in T}$ be a unital field of positive
linear maps $\phi_t :B(H)\to B(K),$ where $K$ is a Hilbert space,
defined on a locally compact Hausdorff space $T$ equipped with a Radon measure $\mu.$
Then for every continuous field of positive operators $(A_t)_{t\in T}$ we have
\begin{equation}\label{DK2}
\int_{T}\phi_t (f(A_t))\,d\mu (t)\leq f(0){\bf 1}+\int_{T}\ff{f(\|
A_t\|)-f(0)}{\| A_t\|}\phi_{t}(A_t)\,d\mu (t).
\end{equation}
Especially, for $f(0)\leq 0,$ the inequality
\begin{equation}\label{DK3}
\int_{T}\phi_t (f(A_t))\,d\mu (t)\leq \int_{T}\ff{f(\| A_t\|)}{\|
A_t\|}\phi_{t}(A_t)\,d\mu (t).
\end{equation}
is valid.
\end{theorem}
\begin{proof} Since $f$ is a convex function, $f(x)\leq
\ff{M-x}{M-m}f(m)+\ff{x-m}{M-m}f(M)$ for every $x\in [m,M]$ where $m\le
M$. Since $\|\cdot\|$ is normalized and unitarily invariant, we have
$0< A_t\leq \|A_t\|{\bf 1}$ and thus
$$f(A_t)\leq \ff{\|A_t\|{\bf 1}-A_t}{\| A_t\|}f(0)+\ff{A_t}{\|
A_t\|}f(\| A_t\| )$$ for every $t\in T$. Applying positive linear
maps and integrating we obtain
\begin{equation}\label{DK4}
\begin{array}{l}
\displaystyle\int_{T}\phi_t(f(A_t))\,d\mu (t)\\[3ex]
\displaystyle\leq f(0)\left[{\bf
1}-\int_{T}\ff{\phi_t(A_t)}{\| A_t\|}\,d\mu
(t)\right]+\int_{T}\ff{f(\| A_t\| )}{\| A_t\|}\phi_t (A_t)\,d\mu(t)
\end{array}
\end{equation}
or
\begin{equation}\label{DK5} \int_{T}\phi_t(f(A_t))\,d\mu (t)\leq
f(0){\bf 1}+\int_{T}\ff{f(\| A_t\| )-f(0)}{\| A_t\|}\phi_t
(A_t)\,d\mu (t).
\end{equation}
Note that since $\int_{T}\ff{\phi_t(A_t)}{\| A_t\|}d\mu (t)\leq
\int_{T}\ff{\| A_t\|\phi_t ({\bf 1})}{\| A_t\|}d\mu (t)={\bf 1}$
we obtain, for $f(0)\leq 0,$ inequality (\ref{DK3}) from (\ref{DK4}) .
\end{proof}

\begin{remark} Setting $T=\{ 1\}$ the inequality (\ref{DK3}) gives
\[
\phi (f(A))\leq \ff{f(\| A\| )}{\| A\|}\phi (A).
\]
Furthermore, setting $\phi=1,$ we get the inequality $f(\| A\|)\geq \| f(A)\|$
obtained in \cite{DK} under the assumption that $f$ is a nonnegative convex function
with $f(0)=0$.

\end{remark}

Related inequalities may be obtained by using subdifferentials. If
$ f\colon\mathbf R\to\mathbf R $ is a convex function and $ [m,M]
$ is a closed bounded real interval, then a subdifferential
function of $ f $ on $[m,M]$ is any function $
k\colon[m,M]\to\mathbf R $ such that
\[
k(x)\in[f'_-(x), f'_+(x)]\qquad x\in(m,M),
\]
where $ f'_- $ and $ f'_+ $ are the one-sided derivatives of $ f $
and $ k(m)=f'_+(m) $ and $ k(M)=f'_-(M). $ Since these functions
are Borel measurable, we may use the Borel functional calculus.
Subdifferential function for concave functions is defined in
analogous way.

\begin{theorem}\label{Sub1} Let $(x_t)_{t\in T}$ be a bounded continuous field
of self-adjoint elements in a unital $C^{*}$-algebra ${\mathcal
A}$ with spectra in $[m,M]$ defined on a locally compact Hausdorff
space $T$ equipped with a Radon measure $ \mu, $ and let
$(\phi_t)_{t\in T}$ be a unital field of positive linear maps
$\phi_t :{\mathcal A}\to {\mathcal B}$ from ${\mathcal A}$ to
another unital $C^{*}$-algebra ${\mathcal B}$. If $f:\R\to\R$ is a
convex function then
\begin{eqnarray}\label{Sub2}
\lefteqn{f(y){\bf 1}+k(y)\left(\int_{T} \phi_t(x_t)d\mu (t)-y{\bf
1}\right)} \nonumber \\
& & \leq \int_{T} \phi_t(f(x_t))\, d\mu (t)  \\
& & \leq f(x){\bf 1}-x\int_{T}\phi_t(k(x_t))\,d\mu (t)+\int_{T}
\phi_t(k(x_t)x_t)\,d\mu (t)\nonumber
\end{eqnarray}
for every $x,y\in [m,M]$, where $k$ is a subdifferential function of $f$ on $ [m,M]. $
In the dual case ($f$ is concave) the opposite inequality holds.
\end{theorem}
\begin{proof} Since $f$ is convex we have $f(x)\geq
f(y)+k(y)(x-y)$ for every $x,y\in [m,M]$. By using the functional
calculus it then follows that $f(x_t)\geq f(y){\bf 1}+k(y)(x_t-y{\bf
1})$ for $t\in T$. Applying the positive linear maps $\phi_t$ and
integrating, LHS of (\ref{Sub2}) follows. The RHS of (\ref{Sub2})
follows similarly by using the functional calculus in the variable $y.$
\end{proof}

Numerous inequalities can be obtained from (\ref{Sub2}). For
example, LHS of (\ref{Sub2}) may be used to obtain an estimation from
below in the sense of Theorem \ref{Rev1}. Namely, the following
theorem holds.
\begin{theorem}\label{Sub3} Let $(x_t)_{t\in T}$ be a bounded continuous field
of self-adjoint elements in a unital $C^{*}$-algebra ${\mathcal
A}$ with spectra in $[m,M]$ defined on a locally compact Hausdorff
space $T$ equipped with a Radon measure $ \mu, $ and let
$(\phi_t)_{t\in T}$ be a unital field of positive linear maps
$\phi_t\colon{\mathcal A}\to {\mathcal B}$ from ${\mathcal A}$ to
another unital $C^{*}$-algebra ${\mathcal B}$. Let $f:\R\to\R$,
$g:[m,M]\to\R$ and $F:U\times V\to\R$ be functions such that
$f\left([m,M]\right)\subset U,$ $g\left([m,M]\right)\subset V,$
$F$ is bounded, $f$ is convex  and $f(y)+k(y)(t-y)\in U$ for every
$y,t\in [m,M],$ where $k$ is a subdifferential function of $f$ on
$ [m,M]. $ If $F$ is operator monotone in the first variable, then
\begin{equation}\label{Sub4}
\begin{array}{l}
\displaystyle F\left[\int_{T}\phi_{t}\left(f(x_t)\right)\,d\mu
(t),g\left(\int_{T}\phi_{t}(x_t)\,d\mu (t)\right)\right]\\[3ex]
\displaystyle\geq
\inf_{m\leq z\leq M}F\left[f(y)+k(y)(z-y),g(z)\right]{\bf 1}
\end{array}
\end{equation}
for every $y\in [m,M]$. In the dual case (when $f$ is concave) the
opposite inequality holds in (\ref{Sub4}) with $\sup$ instead of
$\inf$.
\end{theorem}

Using LHS of (\ref{Sub2}) we can give generalizations of some dual
results from \cite{DK}.

\begin{theorem}\label{DK11} Let $(x_t)_{t\in T}$ be a bounded continuous field
of positive elements in a unital
$C^{*}$-algebra ${\mathcal A}$ defined on
a locally compact Hausdorff space $T$ equipped with a Radon measure
$ \mu, $ and let $(\phi_t)_{t\in T}$ be a unital field of positive
linear maps $\phi_t\colon{\mathcal A}\to {\mathcal B}$ from ${\mathcal
A}$ to another unital $C^{*}$-algebra ${\mathcal B}$ acting on a
finite dimensional Hilbert space $K.$ Let $\|\cdot\|$ be a unitarily
invariant norm on $B(K)$ and let $f\colon [0,\infty )\to \R$ be an increasing function.
\begin{enumerate}
\item If $\|{\bf 1}\|=1$ and $f$ is convex with $f(0)\leq 0$ then
\begin{equation}\label{DK21} f\left(\|\int_{T}\phi_t (x_t)\,d\mu (t)\|\right)\leq
\|\int_{T}\phi_t(f(x_t))\,d\mu (t)\|.\end{equation} \item If
$\int_{T}\phi_{t}(x_t)\,d\mu (t)\leq \|\int_{T}\phi_t(x_t)\,d\mu
(t)\|{\bf 1}$ and $f$ is concave then
\begin{equation}\label{DK31}\int_{T}\phi_t(f(x_t))\,d\mu (t)\leq
f\left(\|\int_{T}\phi_t(x_t)\,d\mu (t)\|\right){\bf
1}.\end{equation}
\end{enumerate}
\end{theorem}
\begin{proof} Since $f(0)\leq 0$ and $f$ is increasing we have $k(y)y-f(y)\geq 0$
and $k(y)\geq 0$. From (\ref{Sub2}) and the triangle inequality we have
\[
k(y)\|\int_{T}\phi_t(x_t)\,d\mu (t)\|\leq\|\int_{T}\phi_{t}(f(x_t))\|+(k(y)y-f(y)).
\]
Now (\ref{DK21}) follows by setting $y=\|\int_{T}\phi_t(x_t)\,d\mu (t)\|$. Inequality
(\ref{DK31}) follows immediately from the assumptions and from the
dual case of LHS in (\ref{Sub2}) by setting $y=\|\int_{T}\phi_t(x_t)\,d\mu(t)\|$.
\end{proof}

Finally, to illustrate how RHS of (\ref{Sub2}) works,  we set
\[
x=\ff{\|\int_{T}\phi_t(k(x_t)x_t)\,d\mu (t)\|}{\|\int_{T}\phi_t(k(x_t))\,d\mu (t)\|}
\]
and obtain a Slater type inequality
\[
\int_{T}\phi_t (f(x_t))\,d\mu(t)\leq
f\left(\ff{\|\int_{T}\phi_t(k(x_t)x_t)\,d\mu
(t)\|}{\|\int_{T}\phi_t(k(x_t))\,d\mu (t)\|}\right){\bf 1}
\]
under the condition
\[
\ff{\int_{T}\phi_t (k(x_t)x_t)\,d\mu (t)}{\|\int_{T}\phi_t
(k(x_t)x_t)\,d\mu (t)\|}\leq \ff{\int_{T}\phi_t (k(x_t))\,d\mu
(t)}{\|\int_{T}\phi_t (k(x_t))\,d\mu (t)\|}\,.
\]

\end{document}